\documentstyle[fleqn,12pt]{article}

\newfont{\ger}{eufm10 scaled \magstep 2}

\def\C{\makebox[.8em][l]{\makebox[.25em][l]{C}\rule{.03em}{1.5ex}}}

\newcommand{\A}{\mbox {{\ger A}} }
\newcommand{\B}{\mbox {{\ger B}} }
\newcommand{\X}{\mbox {{\ger X}} }
\newcommand{\D}{\mbox {${\cal D}$} }
\newcommand{\E}{\mbox {${\cal E}$} }
\newcommand{\U}{\mbox {${\cal U}$} }

\newcommand{\rge}{\rm ran \:}

\newcommand{\Lin}[1] {\mbox {${\cal L}({#1})$}}
\newcommand{\bigfrac}[2]
   {\mbox {$\frac{\displaystyle{{#1}}}{\displaystyle{{#2}}}$  }}

\newcommand{\opmatrix}[4]{\mbox {
    $ \left [ \begin{array} {cc}
                   {#1} & {#2} \\ {#3} & {#4} \end{array} \right ] $}
    }
\newcommand{\col}[2]{\mbox {
    $ \left [ \begin{array} {c} {#1} \\ {#2} \end{array} \right ] $}
    }
\newcommand{\row}[2]{\mbox {
    $  [ \begin{array} {c} {#1} ~ {#2} \end{array}  ] $}
    }
\newcommand{\opmthree}[9]{\mbox {
    $ \left [ \begin{array} {ccc}
                   {#1} & {#2} & {#3} \\
                   {#4} & {#5} & {#6} \\
                   {#7} & {#8} & {#9} \end{array} \right ] $}
    }
\newcommand{\colthree}[3]{\mbox {
    $ \left [ \begin{array} {c} {#1} \\ {#2} \\ {#3} \end{array} \right ] $}
    }
\newcommand{\binom}[2]{\mbox {
    \renewcommand{\arraystretch}{.7}
    $\left(\begin{array}{c}{#1}\\{#2}\end{array}\right)$
    \renewcommand{\arraystretch}{1} }
   }

\newtheorem{thm}{Theorem}
\newtheorem{prop}[thm]{Proposition}
\newtheorem{cor}[thm]{Corollary}
\newtheorem{lem}[thm]{Lemma}

\begin{document}

\title{Unbounded Symmetric Homogeneous Domains in Spaces of Operators}
\author{Lawrence A. Harris}

\date{}
\maketitle

\begin{abstract} We define the domain of a linear fractional transformation in a space of
operators and
show that both the affine automorphisms and the compositions of symmetries
act transitively on these domains.
Further, we show that Liouville's theorem holds
for domains of linear fractional transformations and, with an additional
trace class condition, so does the Riemann removable singularities theorem.
We also show that every biholomorphic
mapping of the operator domain $I < Z^*Z$ is a linear isometry when the
space of operators is a complex Jordan subalgebra of ${\cal L}(H)$
with the removable singularity property
and that every biholomorphic
mapping of the operator domain $I + Z_1^*Z_1 < Z_2^*Z_2$ is a linear
map obtained by multiplication on the left and right by J-unitary and
unitary operators, respectively.
\end{abstract}

{\bf 0. Introduction.}
This paper introduces a large class of finite and infinite dimensional
symmetric affinely homogeneous domains which are not holomorphically
equivalent to any bounded domain.  These domains are subsets of spaces of
operators and include
domains as diverse as a closed complex subspace of the bounded linear
operators from one Hilbert space to another, the identity component of the
group of invertible operators in a $C^*$-algebra and the complement of a
hyperplane in a Hilbert space.
Each of our domains may be characterized as a component of the domain of definition
of some linear fractional transformation which maps a neighborhood of a
point in the component biholomorphically onto an open set in the same space.
Thus we refer to our domains as domains of linear fractional transformations.
We show that at any point of such a domain, there exists a biholomorphic
linear fractional transformation of the domain onto itself which is a symmetry
at the point.  
This linear fractional transformation is a generalization of the
Potapov-Ginzburg transformation.
Moreover, for any two points $Z_0$ and $W_0$ in the domain, there
exists an affine automorphism of the domain of the form
$\phi(Z) = W_0 + A(Z-Z_0)B$ and $\phi$ is a composition of the above
symmetries.
In fact, we exhibit many different formulae for linear fractional
transformations that are automorphisms of domains of linear fractional
transformations.

We show that any bounded holomorphic function on the domain of a linear
fractional transformation is constant and, more generally, that the
Kobayashi pseudometric vanishes identically on such domains.
We also show that if the domain of a linear fractional transformation
satisfies a trace class condition, then a holomorphic function on the
domain which is locally bounded with respect to points outside the
domain has a holomorphic extension to these points.
Further, we give a sufficient condition for two domains
of linear fractional transformations to be affinely equivalent
when the subspace of operators considered is the full space or a $C^*$-algebra.
As expected, our unbounded symmetric homogeneous domains fail to have many
of the properties of bounded symmetric domains.  For example,
Cartan's uniqueness theorem may fail, symmetries may not be unique and
automorphisms may not be a composition of a linear fractional transformation
and a linear map.

Next we turn to a discussion of the characterization of the automorphisms
of non-homogeneous domains
which are operator analogues of the domain $1 + |z_1|^2 < |z_2|^2$ in
$\C^2$ (which is included).
The domains we consider are circular domains in a
space of operators in the sense discussed
previously by the author in \cite{LFT} and are holomorphically equivalent
to the open unit ball of the space with certain singular points
omitted.  We show that with some restrictions,
these points are removable singularities and we deduce that
all biholomorphic mappings of the domains are linear.
In the case where the domain and
range spaces of the operators considered have different dimensions, these
linear maps are given by operator matrices which are the coefficient
matrices of linear fractional transformations of the unit ball of the
space of operators on the domain space.
In the case where the operators considered have the same domain and range,
our domains reduce to analogues of the exterior of the unit disc
in spaces of operators and the linear maps are invertible isometries of these spaces.
A result of Arazy and Solel~\cite{arazy} allows us to treat spaces of operators
which are not necessarily closed under adjoints.

A main lemma is the result that every $J^*$-algebra with identity which is
closed in the weak operator topology has what we call the removable singularity property,
i.e., the singular operators in the open unit ball are removable singularities
for bounded holomorphic functions.
This is rather surprising since these singular operators
may have non-empty interior in the infinite dimensional case.
We also give two types of homogeneous circular domains where the group of linear
automorphisms of the domains acts transitively.

Unbounded homogeneous domains in $\C^n$ have been considered previously by
Penney~\cite{RHD} and Winkelmann~\cite{wink}.  In particular, Penney obtains
a classification of the rationally homogeneous domains in $\C^n$ analogous
to the classification theorem for bounded homogeneous domains given by Vinberg, Gindikin and Piatetskii-Shapiro.
See also Hua~\cite{hua1}.

Our discussion emphasizes the explicit construction of domains and mappings
using ideas from operator theory and functional analysis.  We have considered
bounded symmetric homogeneous domains from this point of view previously
in~\cite{BSD}.
See~\cite{isidro} and~\cite{upmeier} for an exposition of the theory of bounded
symmetric domains in infinite dimensions.

Let $H$ and $K$ be complex Hilbert spaces and let $\Lin{H,K}$ denote the Banach
space of all bounded linear transformations from $H$ to $K$ with the
operator norm.
We write $\Lin{H}$ for $\Lin{H,H}$.
Throughout, $\A$ and $\B$ denote any closed complex subspaces
of $\Lin{H,K}$.
(The reader interested only in the finite dimensional case may take
$H = \C^n$, $K=\C^m$ and identify $\Lin{H,K}$ with the vector space of
all $m \times n$ matrices of complex numbers.)

{\bf 1.  Domains of Linear Fractional Transformations.}
Let $C \in \Lin{K,H}$, $D \in \Lin{H}$ and suppose
$(CZ_0 + D)^{-1}$ exists for some $Z_0 \in \A$.
Put $X_0 = (CZ_0 + D)^{-1} C$.  If
$ZX_0Z \in \A$ for all $Z \in \A$,
we define
\begin{equation}
\label{D:eq}
\D = {\rm Comp}_{Z_0} \{Z \in \A: (CZ+D)^{-1} ~ \mbox{exists} \}
\end{equation}
to be the {\em domain of a linear fractional transformation on} $\A$.
(Here and elsewhere, ${\rm Comp}_{Z_0}$ denotes the connected component
containing $Z_0$.)
It is easy to show that any finite or infinite product of domains of linear
fractional transformations is the domain of a linear fractional
transformation.
To justify this terminology, we observe that by~\cite[Prop. 5]{LFT},
if the linear fractional transformation
\begin{equation}
\label{lft:eq}
T(Z) = (AZ+B)(CZ+D)^{-1}
\end{equation}
maps a domain in $\A$ containing $Z_0$
onto a subset of some $\B$ containing an interior point of $\B$
and if the coefficient
matrix of $T$ is invertible, then  $ZX_0Z \in \A$ for all $Z \in \A$ and
$T$ is a biholomorphic mapping of the domain $\D$ above onto a similar domain in $\B$.
We consider components because a linear fractional transformation may be a
biholomorphic mapping of one component of the open set where it is defined
while mapping the other components outside the space.  (See Example~1 below.)

\begin{thm}
\label{U:thm}
The domain $\D$ of a linear fractional transformation on $\A$
is a symmetric affinely homogeneous domain.  In particular, if $\D$ is
given by~(\ref{D:eq}), then for each $Y \in \D$,
\begin{equation}
\label{U:eq}
U_{Y}(Z) = Y - (Z-Y)(CZ+D)^{-1}(CY+D)
\end{equation}
is a biholomorphic mapping of \D onto itself satisfying $U_{Y}^2 = I$,
$U_{Y}(Y) = Y$ and $DU_{Y}(Y) = -I$.  Moreover, for any point
$W_0 \in \D$, there is an invertible affine linear fractional
transformation $\phi$ of \D onto itself with $\phi(Z_0)=W_0$ and $\phi$ can
be chosen to be a
composition of mappings $U_{Y}$ given above.
\end{thm}

Winkelmann~\cite[Prop. 1]{wink} has proved a related result for the complement of an algebraic
subvariety in $\C^n$.
It is easy to show that $U_Y$ is the only linear fractional transformation
(\ref{lft:eq}) with the mapping property described there which
satisfies $T(Y) = Y$ and $DT(Y) = -I$ on $\A$.
The transformations $U_{Y}$ are used in~\cite{fact}
to obtain factorizations of certain operator matrices.

{\bf Example 0.}
Any $\A$ is the domain of a linear fractional transformation on $\A$ when
$C=0$ and $D=I$ since then $X_0 =0$ for any choice of $Z_0$.  In this case,
$U_Y(Z) = 2Y-Z$.

{\bf Example 1.}  Let $\A$ be a power algebra~\cite{SL}, i.e., $\A$ is a
closed complex subspace of $\Lin{H}$ containing the identity operator $I$
on $H$ and the squares of each of its elements.  Put
\begin{equation}
\label{GI:eq}
G_I(\A) = {\rm Comp}_{I} \{Z \in \A: Z^{-1} ~ \mbox{exists} \}.
\end{equation}
Clearly, $G_I(\A)$ is the domain of a linear fractional transformation on $\A$
since $X_0 = I$ for $Z_0 = I$.
Also, $U_{Y}(Z) = YZ^{-1}Y$.
If $\A$ is a $W^*$-algebra, the argument in~\cite[Cor. 5.30]{BAT} shows that the
invertible operators in $\A$ are connected so it is not necessary to take
components in this case.
However, if $\A$ is the closure of the polynomials in the bilateral
shift~\cite[Example 4.25]{BAT}, then $U_I$ is a biholomorphic mapping of $G_I(\A)$
but $U_I$ does not take any other invertible elements of $\A$ (such as the
bilateral shift) to $\A$.

{\bf Example 2.}
Let $\A$ be a closed complex subspace of $\Lin{H}$ containing a projection
$E$ such that $ZEZ \in \A$ whenever $Z\in\A$.  Then
\[
\D = {\rm Comp}_{E} \{Z \in \A: (EZ+I-E)^{-1} ~ \mbox{exists} \}
\]
is the domain of a linear fractional transformation on $\A$ since $X_0 = E$
for $Z_0 = E$. The previous examples are the cases $E=0$ and $E=I$.
Moreover,
\[
U_E(Z) = [(E-I)Z + E](EZ+I-E)^{-1}
\]
is a variant form of the Potapov-Ginzburg transformation~\cite[\S 3]{alpay}.
(Another variant is called the Redheffer transform in~\cite[p. 269]{BGR}.
An early reference is~\cite[p. 240]{hua2}.)
A basic property is that $U_E$ is a biholomorphic mapping of
$\D_1 \cap \D$ onto $\D_2 \cap \D$, where
\[
\D_1  =  \{Z \in \A:~Z^*JZ < J \}, ~~~~
\D_2  =  \{Z \in \A:~||Z|| < 1 \},
\]
and $J=I-2E$.
This follows directly from~\cite[Lemma 7]{LFT} since the coefficient matrix
$M$ of $U_E$ satisfies $J_2 = (M^{-1})^*J_1M^{-1}$, where
\[
J_1 = \opmatrix{J}{0}{0}{-J},~~~~~~J_2 = \opmatrix{I}{0}{0}{-I}.
\]

{\bf Example 3.}
If $\rge C$ is closed or if $\rge D \subseteq \rge C$, then
\[
 \D = \{Z \in \Lin{H,K}:~(CZ+D)^{-1}~{\rm exists} \}
\]
is the domain of a linear fractional transformation on $\A = \Lin{H,K}$
when $\D$ is not empty
since $\A$ is obviously closed under the required products and since
$\D$ is connected by Proposition~\ref{connect:prop} below.

{\bf Example 4.}
Let $c \in H$ with $c \neq 0$ and $d \in \C$.  Then
\[
\D = \{z \in H:  (z,c) \neq -d \}
\]
is the complement of a hyperplane in $H$ and
the domain of a linear fractional transformation on $H$
with
\[
 U_y(z) = y - \frac{(y,c)+d}{(z,c)+d}(z-y).
\]
(Identify $H$ with $\Lin{\C,H}$ and apply Example~3 with $C = c^*$.)

{\bf Example 5.}
Let $\A$ contain all the rank one operators in $\Lin{H,K}$ and let
$x \in H$, $y \in K$ be unit vectors.  Then if $d \neq 0$,
\[
\D = \{Z \in \A:  (Zx,y) \neq -d \}
\]
is the domain of a linear fractional transformation on $\A$.  Indeed,
let $C = xy^*$, $D = dI$ and note that $CZ + D = d I + x(y^*Z)$ is
invertible if and only if $d + (y^*Z)x$ is invertible.
Also $X_0 = xy^*$ for $Z_0 = (1-d)yx^*$ and
\[
 U_Y(Z) = Y - \frac{Z-Y}{d}\left [I - \frac{xy^*Z}{d+(Zx,y)}\right ](dI+xy^*Y).
\]
Since $C$ is compact, $\D$ is connected by Proposition~\ref{connect:prop}
below.

{\bf Example 6.}
Let $H$ be a Hilbert space with conjugation $z \rightarrow \bar{z}$ and
let $c \in H$ with $c \neq 0$.  Then each of the domains
\begin{eqnarray*}
 \D_1 &=& \{z \in H:~1 + 2(z,c) + (z,\bar{z})(\bar{c},c) \neq 0 \},\\
 \D_2 &=& \{z \in H:~(z,\bar{z}) \neq 0 \}
\end{eqnarray*}
is linearly equivalent to the domain of a linear fractional transformation.
To see this, observe that by~\cite[\S2]{BSD} there exists an invertible linear
map $z \rightarrow A_z$ of $H$ onto a Cartan factor $\A$ of type IV.  For $\D_1$,
take $C=A_c^*$, $Z=A_z$ and $D=I$, and note that $CZ+D$ is invertible if and
only if $z\in\D_1$ by~\cite[(9)]{BSD}.
Also, $X_0 = A_c^*$ for $Z_0 = 0$ so $ZX_0Z \in \A$ for all $Z\in\A$ since
$\A$ is a $J^*$-algebra.
Moreover, $\D_1$ is connected by Proposition~\ref{connect:prop} below
since $X_0Z$ satisfies a quadratic polynomial equation.
A similar argument with $C=I$ and $D=0$ establishes the result
for $\D_2$, and in that case,
\[
 U_y(z) = \frac{2(z,\bar{y})y -(y,\bar{y})z}{(z,\bar{z})}.
\]

{\bf Proof of Theorem 1.}
Given $Y \in \D$, put $X = (CY + D)^{-1}C$ and note the identities
\begin{eqnarray}
\label{CU+D:eq}
CU_{Y}(Z)+D &=& (CY+D)(CZ+D)^{-1}(CY+D),\\
\label{X:eq}
(CY + D)^{-1}(CZ + D) &=& I + X(Z - Y).
\end{eqnarray}
The coefficient matrix of $U_{Y}$ is
\[
M = \opmatrix{-(I-YX)}{2Y-YXY}{X}{I-XY}
\]
by~(\ref{X:eq}).
Clearly $M^2 = I$ so M is invertible and $U_{Y}^2 = I$.  Also,
$U_{Y}(Y) = Y$.  Hence by hypothesis and \cite[Prop. 5]{LFT},
$U_{Y}$ is a biholomorphic mapping of $\D$ onto itself when $Y = Z_0$
and therefore $ZXZ \in \A$ whenever $Z \in \A$.  Thus \cite[Prop. 5]{LFT}
applies again to show that $U_{Y}$ is a biholomorphic mapping of $\D$ onto
itself for arbitrary $Y \in \D$.  Further, $DU_{Y}(Y) = -I$ since
$DU_{Y}(Y)Z = \bigfrac{d}{dt} U_{Y}(Y + tZ) |_{t=0}$ and
$U_{Y}(Y + tZ) = Y -tZ(I+tXZ)^{-1}$.

Next we show that if $Z, W \in \D$ satisfy
$||(CZ+D)^{-1}C(W-Z)|| < 1$, then there is a $Y \in \D$ with
$U_{Y}(Z) = W$.  Put $X = (CZ+D)^{-1}C$ and $R=W-Z$.
By the holomorphic functional calculus~\cite[\S5.2]{HP}, the binomial series defines a
$Q \in \Lin{H}$ such that $Q^2 = I + XR$, $\sigma(Q)$ does not contain $-1$,
and $(I+Q)^{-1}$ is a limit of polynomials in $XR$.
Take $Y = Z + R(I+Q)^{-1}$.  Then $Y \in \A$ by \cite[Lemma 6a]{LFT}.
Also, since $XR = Q^2 -I$, we have that
\begin{eqnarray*}
U_{Y}(Z) - Z & = & Y-Z + (Y-Z)[I+X(Y-Z)]  \\
               & = & R(I+Q)^{-1}[2(I+Q)+XR](I+Q)^{-1}  \\
               & = & R
\end{eqnarray*}
so $U_{Y}(Z)=W$.
Moreover, $(CY+D)^{-1}$ exists by this and~(\ref{CU+D:eq}).
To show that $Y \in \D$, let $0\leq t \leq 1$.  The above argument
applies to the points $Z$ and $W_t = (1-t)Z+tW$ and the corresponding
$Q_t$ is a continuous function of $t$ in $[0,1]$.  Therefore, there is a curve
$\gamma$ in $\A$ connecting $Z$ to $Y$ such that $C\gamma(t)+D$ is invertible
for all $0\leq t \leq 1$.

Now given $W_0 \in \D$, there is a curve $\gamma$ in $\D$ connecting
$Z_0$ to $W_0$.  By compactness, there is a number $M$ such that
$||(C\gamma(t) + D)^{-1}C|| \leq M$ for all $0\leq t \leq 1$ and there is
a $\delta > 0$ such that $||\gamma(s) - \gamma(t)|| < 1/M$ whenever
$|s-t| \leq \delta $ and $0 \leq s,t \leq 1$.  Choose $n > 1/\delta$
and put $Z_k = \gamma(k/n)$ for $0 \leq k \leq n$.  Then for $1 \leq k \leq n$,
\begin{equation}
\label{pts:eq}
||(CZ_{k-1}+D)^{-1}C(Z_k-Z_{k-1})|| < 1,
\end{equation}
so by what we have shown there exists a $U_k$ with $U_k(Z_{k-1}) = Z_k$.
Thus $U = U_1 \cdots U_n$ satisfies
$U(Z_0) = W_0$.

To complete the proof, it suffices to show that 
$U_W \circ U_Y$ is an affine linear fractional transformation for any
$W,Y \in \D$ since then we can choose $n$ to be even in the above argument
and take $\phi=U$.
Put $X = (CY+D)^{-1}C$ and set
\[
  Q = [W -U_Y(Z)][I+X(Z-Y)].
\]
Then
  $Q  = W-Y + [I + (W-Y)X](Z-Y)$
by~(\ref{X:eq}).
Hence by~(\ref{CU+D:eq}) and~(\ref{X:eq}),
\begin{eqnarray*}
U_W \circ U_Y(Z) &=& W  + Q(CY+D)^{-1}(CW+D)  \\
                 &=& U_W(Y) + [I+(W-Y)X](Z-Y)[I+X(W-Y)],
\end{eqnarray*}
as required.    \hfill $\Box$

There are many biholomorphic mappings of the domain of a linear fractional
transformation on $\A$ besides~(\ref{U:eq}).  For example, if $W_0 \in \D$
and $||X_0(W_0-Z_0)|| < 1$, then
\[
\phi(Z) = W_0 +[I+(W_0-Z_0)X_0]^{1/2}(Z-Z_0)[I+X_0(W_0-Z_0)]^{1/2}
\]
is an invertible affine linear fractional transformation of $\D$ onto itself
with $\phi(Z_0) = W_0$, where the square roots are defined by the binomial
series.
This follows from~\cite[Lemma 6c]{LFT}, (\ref{X:eq}) and the identity
\[
 I+X_0[\phi(Z)-Z_0] = R^{1/2}[I+X_0(Z-Z_0)]R^{1/2},
\]
where $R = I+X_0(W_0-Z_0)$.
Moreover, the above proof shows that finite compositions of these
mappings act transitively on $\D$.
Note that $\phi = U_W \circ U_Y$ when $W_0=U_W(Y)$ and $Z_0 = Y$.
If $\A = \Lin{H,K}$ or if $\A$ is a $C^*$-algebra with $C,D \in \A$, a simpler formula for a transitive set of
affine mappings of $\D$ is given by~(\ref{Raff:eq}) below with $R=I$.

There is a general class of involutory biholomorphic mappings of $\D$ which
contains the transformations~(\ref{U:eq}).  Specifically, if $W_0 \in \D$
and $||X_0(W_0-Z_0)|| < 1$, define
\[
V_{W_0}(Z) = Z_0 - T_{W_0-Z_0}(Z-Z_0),
\]
where
\[
T_A(Z) = (I+AX_0)^{-1/2}(Z-A)(I+X_0Z)^{-1}(I+X_0A)^{1/2}.
\]
(Compare~\cite[p. 146]{OSD}.)
Note the identities $V_{W_0} = \phi \circ U_{Z_0}$ and
$U_{Z_0} = V_{W_0} \circ \phi$.  It follows that
$V_{W_0}$ is a biholomorphic mapping of $\D$ onto itself with
$V_{W_0}^2 = I$ and $V_{W_0}(Z_0) = W_0$ and that
the composition $V_{W_1} \circ V_{W_0}$ of any two such mappings
is affine.
In particular, $V_{Z_0} = U_{Z_0}$.

The domain of a linear fractional transformation may have biholomorphic
mappings which cannot be expressed as the composition of a linear fractional
transformation and an invertible linear map.  For example, let
$\A = \Lin{\C,\C^2}$, $C=\row{0}{1}$ and $D = 0$.  Then the domain $\D$
of~(\ref{D:eq}) is 
\[
 \D = \{\col{z_1}{z_2}\in\A:~z_2 \neq 0\}
\]
and
\[
 h\col{z_1}{z_2} = \col{z_1+g(z_2)}{z_2}
\]
is a biholomorphic mapping of $\D$ for any holomorphic function
$g:\C\setminus\{0\}\rightarrow \C$.
Note that Cartan's uniqueness theorem fails for $\D$ since if
$g(z) = (z-1)^2$ and $Z_0 = \col{0}{1}$, then $h(Z_0) = Z_0$ and
$Dh(Z_0) = I$, but $h \neq I$ on $\D$.
Moreover, $U_{Z_0}$ and $h^{-1}\circ U_{Z_0}\circ h$ are distinct
symmetries of $\D$ at $Z_0$.

Let $\X$ be any Banach space.
\begin{thm}
\label{liou:thm}
If $\D$ is the domain of a linear fractional transformation on $\A$ and if
$h:\D \rightarrow \X$ is a bounded holomorphic function, then $h$ is
constant.  In particular, $\D$ is not holomorphically equivalent to a
bounded domain.
\end{thm}

{\bf Proof.}
Let $Z \in \A$ with $||X_0(Z-Z_0)|| < 1$.  We will construct an entire function
$f:\C \rightarrow \D$ with $f(0) = Z_0$ and $f(1) = Z$.  Then $g(\lambda) =
h(f(\lambda))$ is a bounded entire function so $h(Z) = h(Z_0)$ since $g$ is
constant \cite[Th. 3.13.2]{HP}.  Hence $h$ is constant in a
neighborhood of $Z_0$ so $h$ is constant on $\D$ by the identity
theorem~\cite[Th. 3.16.4]{HP}.

Recall~\cite[\S6.3]{AFT} that
$\limsup_{n \rightarrow \infty} \left |\binom{\lambda}{n}\right| ^{1/n} \leq 1$
and that the binomial series
$b_{\lambda}(z) = \sum_{n=0}^{\infty} \binom{\lambda}{n} z^n$
satisfies $b_{\lambda}(z) b_{-\lambda}(z) = 1$ for all $\lambda \in \C$ and
$|z| < 1$.  Put $W = X_0(Z-Z_0)$.  By the holomorphic functional calculus,
$b_{\lambda}(W)b_{-\lambda}(W) = I$ and hence $b_{\lambda}(W)^{-1}$ exists
for all $\lambda \in \C$.  Define
$f(\lambda) = Z_0 + (Z-Z_0) \sum_{n=1}^{\infty} \binom{\lambda}{n} W^{n-1}$.
Then $f$ is entire and $f(\lambda) \in \A$ for all $\lambda \in \C$
by~\cite[Lemma 6a]{LFT}.  Moreover, if $\lambda \in \C$, then
\[
(CZ_0+D)^{-1}[Cf(\lambda)+D] = I + X_0[f(\lambda)-Z_0] = b_{\lambda}(W),
\]
so $f(\lambda) \in \D$.   \hfill $\Box$

Note that it follows from the properties of $f(\lambda)$ given above and~(\ref{pts:eq}) that
the pseudometric $\rho$ assigned to $\D$ by any Schwarz-Pick System~\cite{SPS}
satisfies $\rho \equiv 0$.  Thus, in particular, the Caratheodory and 
Kobayashi pseudometrics for $\D$ vanish identically.
(Compare~\cite[Cor. 1]{wink}.)

The next result gives a class of domains where all points outside the domain
are removable singularities for locally bounded holomorphic functions.

\begin{thm}
\label{RS:thm}
Suppose the set $\D = \{Z\in \A:~(CZ+D)^{-1}~ {\rm exists}\}$
is non-empty and $C$ is a trace class operator.
Then $\A \setminus \D$ is an analytic set.  Moreover, if $\E$ is
a domain in $\A$ and if $h:\D \cap \E \rightarrow \X$ is a holomorphic
function which is locally bounded in the sense that each operator in
$\E \setminus \D$ has a
neighborhood $\U$ such that $h$ is bounded on $\U \cap \D$,
then $h$ extends to a holomorphic function $\hat{h}:\E \rightarrow \X$.
\end{thm}

It is easy to show that the conclusions of Theorem~\ref{liou:thm} hold for any
$\D$ satisfying the hypotheses of Theorem~\ref{RS:thm}.
See~\cite{ramis} for the definition and properties of analytic sets in
Banach spaces.
Without some restriction on $C$ or $\A$,
the complement of $\D$ may have non-empty interior and thus
is not an analytic set.  For example, if there is a $Z_1 \in \A$ such that
$CZ_1+D$ has only a one-sided inverse $Q$, then an elementary argument shows
that $CZ+D$ is not invertible for
any $Z \in \A$ with $||Z-Z_1|| < 1/r$, where $r = ||C|| \: ||Q||$.

{\bf Proof.}
By~(\ref{X:eq}) and a translation, we may suppose that $D=I$.
It suffices to show that there is a holomorphic function $f:\A \rightarrow \C$
satisfying
\begin{equation}
\label{A:eq}
 \D = \{Z\in\A:~f(Z) \neq 0\}
\end{equation}
since then $\A \setminus \D$ is an
analytic set and the remainder of the theorem follows
from the extension of the Riemann removable singularities theorem
given in~\cite[Th. II.1.1.5]{ramis}.
It is easy to construct such an $f$ when $C$ has finite rank.  Indeed,
put $H_1 = \rge C$, $K_1 = (\ker C)^{\bot}$,
$H_2 = H_1^{\bot}$ and $K_2 = K_1^{\bot}$.
Since these spaces are closed, the operators in $\Lin{H,K}$
and $\Lin{K,H}$ can be written as associated $2\times 2$ operator matrices.
In particular, if $Z \in \Lin{H,K}$,
\begin{equation}
\label{dec:eq}
Z = \opmatrix{Z_1}{Z_2}{Z_3}{Z_4}, ~~~ C = \opmatrix{Y}{0}{0}{0},
\end{equation}
where $Z_1 \in \Lin{H_1,K_1}$, $Z_2 \in \Lin{H_2,K_1}$,
      $Z_3 \in \Lin{H_1,K_2}$, $Z_4 \in \Lin{H_2,K_2}$
and $Y \in \Lin{K_1,H_1}$.  Hence
\begin{equation}
\label{i+xz:eq}
I + CZ = \opmatrix{I + YZ_1}{YZ_2}{0}{I},
\end{equation}
where the successive appearances of $I$ denote the identity operator
on $H$, $H_1$ and $H_2$, respectively.
Thus $Z \in \D$ if and only if $(I+YZ_1)^{-1}$ exists.
Since $H_1$ is finite dimensional, 
we may take $f(Z) = \det(I+YZ_1)$ to obtain~(\ref{A:eq}).

If $C$ is a trace class operator, there exists a sequence $\{F_n\}$
of finite rank operators in $\Lin{H,K}$ with $||F_n - C||_1 \rightarrow 0$,
where $||~~||_1$ is the trace class norm~\cite[Th. VI.4.1]{GGK}.
Let $f_n$ be the function defined above where $C$ is replaced by $F_n$.
Then by~\cite[p.115-117]{GGK}, the sequence $\{f_n\}$ is locally uniformly
bounded on $\A$ and converges pointwise to a function $f$ on $\A$ which is
denoted by $f(Z) = \det(I+CZ)$.
Hence $f$ is holomorphic on $\A$ by an extension of Vitali's
theorem~\cite[Th. 3.18.1]{HP} and (\ref{A:eq}) holds by~\cite[VII.3(13)]{GGK}.
\hfill $\Box$

In some cases it is not necessary to take components in the definition
of the domain of a linear fractional transformation.

\begin{prop}
\label{connect:prop}
The set
\begin{equation}
\label{GD:eq}
 \D = \{Z\in \A:~(CZ+D)^{-1}~ {\rm exists}\}
\end{equation}
is connected if at least one of the following conditions holds:
\begin{enumerate}
\item[a)] $C\A$ contains only compact operators,
\item[b)] there is some $Z_0 \in \D$ such that each operator in $X_0\A$ satisfies a non-trivial polynomial equation,
\item[c)] $\A = \Lin{H,K}$ and $\rge C$ is closed,
\item[d)] $\A = \Lin{H,K}$ and $\rge D \subseteq \rge C$.
\end{enumerate}
\end{prop}

Here ran denotes the range of operators.
It is not difficult to deduce from (\ref{X:eq}) and \cite[Lemma 18]{LFT} that
the domain $\D$ above is a circular
domain in the sense of~\cite{LFT} when (a) or (b) holds.
(Take $J = - \row{C}{D}^* \row{C}{D}$.)

{\bf Proof.}
Suppose that (a) holds and that $\D$ is not empty.
As in the proof of Theorem~\ref{RS:thm}, we may suppose that $D=I$.
Let $Z \in \D$ and put
$\Omega = \{\mu \in \C:~\mu Z \in \D\}.$
Clearly $\Omega$ contains $0$ and $1$, and the set $\C\setminus\Omega$
is discrete by hypothesis.
Hence there is a curve
$\omega$ in $\Omega$ connecting $0$ and $1$ so $\gamma = \omega Z$
is a curve in $\D$ connecting $0$ and $Z$.  Therefore, $\D$ is connected.
A similar argument shows that $\D$ is connected when (b) holds.

Suppose that (c) holds.
As in the proof of Theorem~\ref{RS:thm}, we may suppose that $D=I$.
Moreover, the decompositions~(\ref{dec:eq})
hold, where $Y\in\Lin{K_1,H_1}$ is invertible and $H_1$ is a closed
subspace of $H$.
Let $Z \in \D$.  Then $(I+YZ_1)^{-1}$ exists so by the
connectedness~\cite[Cor. 5.30]{BAT} of the set $G$ of invertible
elements of $\Lin{H_1}$, there exists a curve $\omega$ in
$G$ which connects $I$ to $I+YZ_1$.  Define
\[
\gamma(t) = \opmatrix{Y^{-1}\omega(t)-Y^{-1}}{tZ_2}{tZ_3}{tZ_4},
~~~ 0\leq t \leq 1.
\]
Then $I+C\gamma(t)$ is invertible for $0\leq t \leq 1$ by~(\ref{i+xz:eq})
so $\gamma(t)$ is a curve in $\D$ which connects $0$ to $Z$.  Thus $\D$ is
connected.

Suppose that (d) holds and that $Z_0 \in \D$.  By~\cite{RGE}, there is a
$Z_1 \in \Lin{H,K}$ with $CZ_1 = D$ so $C(Z_1+Z_0)$ is invertible.
Then $\rge C = H$ and hence $\D$ is connected by part (c).   \hfill $\Box$

The following gives a condition under which two domains of the
form~(\ref{GD:eq}) are affinely equivalent.

\begin{prop}
\label{a_equiv:prop}
Let $C_1, C_2 \in \Lin{K,H}$, $D_1, D_2 \in \Lin{H}$ and put
\begin{eqnarray*}
  \D_1 &=& \{Z \in \Lin{H,K}: ~(C_1Z+D_1)^{-1}~ {\rm exists} \},\\
  \D_2 &=& \{Z \in \Lin{H,K}: ~(C_2Z+D_2)^{-1}~ {\rm exists} \}.
\end{eqnarray*}
Suppose there exists an invertible $R \in \Lin{K}$ with $C_2 = C_1R$.
Then for each $Z_1 \in \D_1$ and $Z_2 \in \D_2$
there is an invertible affine linear fractional
transformation $\phi$ of $\D_1$ onto $\D_2$ with $\phi(Z_1) = Z_2$.
\end{prop}

{\bf Proof.}
Observe that
\begin{equation}
\label{Raff:eq}
\phi(Z) = Z_2 + R^{-1}(Z-Z_1)(C_1Z_1+D_1)^{-1}(C_2Z_2+D_2)
\end{equation}
is the required affine mapping since
\[
C_2\phi(Z)+D_2 = (C_1Z+D_1)(C_1Z_1+D_1)^{-1}(C_2Z_2+D_2).
\]
\vspace{-.3in}
\hfill $\Box$

Note that Proposition~\ref{a_equiv:prop} and its proof also hold when all
the spaces of operators mentioned are replaced by the same $C^*$-algebra.

{\bf 2.  Linear Automorphisms of Domains.}
The main results of this section are that certain unbounded circular domains
(in the sense of \cite{LFT}) have only linear automorphisms.
In the finite dimensional case, the domains we consider are incomplete matrix
Reinhardt domains and extensions of these to non-square matrices.  (See~\cite{MRD}.)
Throughout, we let $\A_0$ denote the open unit ball of $\A$ in the operator
norm.  Thus
\[
 \A_0 = \{Z\in\A:~||Z|| < 1\}.
\]
\begin{thm}
\label{LU:th}
Suppose $H$ is finite dimensional and $K \neq \{0\}$.
Let $\A = \Lin{H,K \times H}$
and put
\[
\D = \{ Z \in \A: I + Z_1^*Z_1 < Z_2^*Z_2 \},
\]
where $Z = \col{Z_1}{Z_2}$ is the decomposition of $Z$ with
$Z_1 \in \Lin{H,K}$ and $Z_2 \in \Lin{H}$.
Then $h$ is a biholomorphic map of $\D$ onto itself if and only $h(Z) = LZU$, where $L$
is an invertible linear map in $\Lin{K \times H}$ satisfying $L^*JL = J$,
$J = \opmatrix{I}{0}{0}{-I}$ and $U \in \Lin{H}$ is unitary.
\end{thm}

\begin{cor}
\label{hT:cor}
Suppose $K \neq \{0\}$.  Then
$h\col{z}{\lambda} = \opmatrix{A}{b}{c^*}{d} \col{z}{\lambda}$
is a biholomorphic mapping of
\[
\D = \{ \col{z}{\lambda}\in K \times \C:~1+||z||^2 < |\lambda|^2 \}
\]
if and only if
$T(z) = \bigfrac{Az+b}{(z,c)+d}$
is a biholomorphic mapping of $K_0 = \{z\in K:~||z|| < 1\}$, where the
coefficients of $T$ have been multiplied by a positive number so that
$|d|^2 - ||b||^2 = 1$.
Moreover, there are no other
biholomorphic mappings of $\D$ or $K_0$.
\end{cor}

Corollary~\ref{hT:cor} follows from Theorem~\ref{LU:th}, the
form of biholomorphic mappings of $K_0$, and the uniqueness (up to a complex
scalar) of the coefficients of linear fractional transformations.
(See (\ref{LFT:eq}) below and extend~\cite[Th. 3]{pota}.)
Note that Theorem~\ref{LU:th} is no longer true when $K = \{0\}$ (i.e., $Z_1$ does
not appear) since $h(Z) = Z^t$ is a biholomorphic mapping of $\D$ which is
not of the specified form.
We consider the case $K=\{0\}$ next.

Let $\A$ be a power algebra.
We say that $\A$ has the {\em removable singularity property} if every
bounded holomorphic function $h:G_I(\A)\cap\A_0 \rightarrow \X$ extends to
a holomorphic function on $\A_0$ for any Banach space $\X$.
It follows from~\cite{bogdan} that there is no loss of generality if one
takes $\X = \C$.  Thus the classical Riemann removable singularities theorem
for several complex variables shows that $\A$ has the removable singularity
property whenever $\dim \A < \infty$.
Recall~\cite{GCA} that a $J^*$-algebra with identity may be characterized
as a power algebra which contains the adjoints of each of its elements.

\begin{prop}
\label{RSJ:prop}
Every $J^*$-algebra with identity which is closed in the weak operator topology has the
removable singularity property.
\end{prop}

\begin{thm}
\label{ext:th}
Let $\A$ be a power algebra with the removable singularity
property and put
\[
\D = \{Z\in G_I(\A):~I<Z^*Z \}.
\]
Then $h$ is a biholomorphic mapping of $\D$ onto itself if and only if $h = L$,
where $L$ is a linear isometry of $\A$ onto itself with $L(I) \in G_I(\A)$.
\end{thm}

Thus by Proposition~\ref{RSJ:prop}, the conclusions of Theorem~\ref{ext:th}
hold, for example, for
$\Lin{H}$, any $W^*$-algebra, any finite rank $J^*$-algebra having a unitary
element, any Cartan factor of type II or IV and any Cartan factor of type III
where the dimension of the underlying Hilbert space is even or infinite.
Note that the condition $I <Z^*Z$ in the definition of $\D$ above can also be
written as $I <ZZ^*$ since for $Z \in G_I(\A)$, the identity
 $||(Z^{-1})^*|| = ||Z^{-1}||$
holds
and  $I <Z^*Z$ if and only if $||Z^{-1}|| < 1$.
The domains $\D$ of Theorems~\ref{LU:th} and \ref{ext:th} are not
homogeneous.  To see this when $\D$ is as in Theorem~\ref{LU:th},
observe that $Z_r=\col{0}{rI}$ is in $\D$ for $r>1$.  By Theorem~\ref{LU:th},
if $h$ is a biholomorphic mapping of $\D$ with $h(Z_r)=Z_s$, then
$I + Z_s^*JZ_s = U^*(I + Z_r^*JZ_r)U$
so $s = r$.  The case of Theorem~\ref{ext:th} is similar.

{\bf Proof of Theorem~\ref{LU:th}.}
Suppose $h$ has the given form $h(Z) = LZU$.  Then
$h(Z)^*Jh(Z) = U^*Z^*JZU$ and
\[
\D = \{Z\in\A:~I+Z^*JZ < 0\},
\]
so $h(\D) \subseteq \D$.
Since $h^{-1}(Z) = L^{-1}ZU^{-1}$ has
the same form as $h$, it follows that $h$ is a biholomorphic mapping of $\D$.

Now suppose $h$ is a biholomorphic mapping of $\D$ and put
$\E = \{Z\in\A:~Z_2^{-1}~{\rm exists} \}$.
Clearly $\D \subseteq \E$
since $\dim H <\infty$ and note that $Z \in \A_0$ if and only if $Z\in\A$
and $Z^*Z < I$.  Hence $T\col{Z_1}{Z_2} = \col{Z_1Z_2^{-1}}{Z_2^{-1}}$
is a biholomorphic mapping of $\D$ onto $\A_0 \cap \E$ with $T^{-1} = T$.
Therefore, 
$g=T\circ h\circ T$ is a biholomorphic mapping of $\A_0 \cap \E$.
By Theorem~\ref{RS:thm} with $C = \row{0}{I}$,
$g$ extends to a holomorphic mapping
$\hat{g}:\A_0 \rightarrow \A$.
To show that $\hat{g}(\A_0) \subseteq \A_0$, suppose this is false.
Then by the maximum principle \cite[Th. 3.18.4]{HP},
there is a $Z_0 \in \A_0$ with $||\hat{g}(Z_0)|| > 1$.  Let $u, v\in H$ be
unit vectors with $(\hat{g}(Z_0)u,v) > 1$ and define
\[
f(Z) = \frac{1}{ \left ( [\hat{g}(Z) -\hat{g}(Z_0)]u,v \right )}.
\]
Then $f$ is holomorphic on $\A_0 \cap \E$ but $f$ does not have a holomorphic extension
to $\A_0$, contradicting Theorem~\ref{RS:thm}.
Similarly, $g^{-1}$ extends
to a holomorphic map of $\A_0$ into $\A_0$ which is $\hat{g}^{-1}$ by
the identity theorem~\cite[Th. 3.16.4]{HP}.  Hence $\hat{g}$ is a biholomorphic
mapping of $\A_0$ which takes $\A_0 \setminus \E$ onto itself.
(Compare~\cite[Lemma 2.4.4]{vigue}.)

By~\cite[Th. 5.3]{franz}, we may write
\begin{equation}
\label{LFT:eq}
 \hat{g}(Z) = (AZ+B)(CZ+D)^{-1},~~~~~M = \opmatrix{A}{B}{C}{D},
\end{equation}
where $M$ is an invertible operator in 
$\Lin{(K\times H)\times H}$
satisfying $M^*J_1M = J_1$ and $J_1 = \opmatrix{I_{K\times H}}{0}{0}{-I_H}$.
Here $A \in \Lin{K\times H}$, $B \in \A$, $C\in\Lin{K\times H,H}$ and
$D\in\Lin{H}$.  Let
\[
A = \opmatrix{A_1}{A_2}{A_3}{A_4},~~~~
B = \col{B_1}{B_2},~~~~
C = \row{C_1}{C_2}
\]
be the corresponding decompositions.
It suffices to show that $A_2 = 0$, $A_3 = 0$, $B_2 = 0$, $C_2 = 0$ and
$A_4^{-1}$ exists; for then the identity $h = T\circ g\circ T$ gives
$h(Z) = LZU$, where $L = \opmatrix{A_1}{B_1}{C_1}{D}$ and $U=A_4^{-1}$,
and it follows from $M^*J_1M = J_1$ that $L$ and $U$ are as asserted.

Pick an $r$ with $0 < r < 1$ and fix $Z_1 \in \Lin{H,K}$ with
$||Z_1|| \leq \sqrt{1-r^2}$.  Define $\phi(Z_2) = A_3Z_1 + A_4Z_2 + B_2$
for $Z_2 \in \Lin{H}$.
By the mapping properties of $g$, if $||Z_2|| < r$, then $\phi(Z_2)$ is
invertible whenever $Z_2$ is invertible and $\phi(Z_2)$ is singular
whenever $Z_2$ is singular.
Since the determinant is analytic, it follows from the identity theorem
that $\phi(Z_2)$ is singular for all singular $Z_2$.

If $A_4$ is not invertible, there is a unit vector $u \in H$ with
$A_4u = 0$.  Let $Z_2$ be invertible with $||Z_2|| < r$ and let
$Z_2' = (I-uu^*)Z_2$.
Then $Z_2'$ is singular since $u^*Z_2' = 0$, and $\phi(Z_2')=\phi(Z_2)$,
a contradiction.  Thus $A_4$ is invertible.
Put $B=A_4^{-1}(A_3Z_1+B_2)$.  Then $Z_2+B$ is singular whenever $Z_2$
is singular.  To show that $B=0$, let $u\in H$ be a unit vector and
take $Z_2 = (I-B)(I-uu^*)$.  Then $Z_2$ is singular since $Z_2u=0$ so
$Z_2 + B = I + (Bu - u)u^*$ is singular.  Hence,
$(Bu,u) = 1 + u^*(Bu-u) = 0$.

We have shown that $A_4^{-1}$ exists and that $A_3Z_1 + B_2 = 0$ for
$||Z_1|| < \sqrt{1-r^2}$, so $A_3 = 0$ and $B_2 = 0$.  Now the coefficient
matrix of $\hat{g}^{-1}$ is $M^{-1} = \opmatrix{A^*}{-C^*}{-B^*}{D^*}$ and hence
the above argument applied to $\hat{g}^{-1}$ shows that $A_2 = 0$ and
$C_2=0$, as required.    \hfill $\Box$

We deduce Theorem~\ref{ext:th} from the following Lemma.

\begin{lem}
\label{L:lem}
Let $\A$ be a power algebra and
let $L$ be a linear isometry of $\A$ onto itself.  Then $U=L(I)$ is
unitary, $U^* \in \A$, $L(G_I(\A))$ is a component of the set $G(\A)$ of
invertible operators in $\A$ and
\begin{equation}
\label{Linv:eq}
L(Z^{-1}) = UL(Z)^{-1}U
\end{equation}
for all $Z\in G_I(\A)$.
\end{lem}

{\bf Proof. }
By~\cite[Cor. 2.8]{arazy}, $U^* \in \A$ so by~\cite[Th. 2]{SL}, $U$ is
unitary and $L(Z) = U\rho(Z)$, where $\rho:\A \rightarrow U^*\A$ is an
invertible linear isometry satisfying $\rho(I) = I$ and
$\rho(Z^2) = \rho(Z)^2$ for all $Z \in \A$.  Hence
\begin{equation}
\label{ptriple:eq}
\rho(ZWZ) = \rho(Z)\rho(W)\rho(Z)
\end{equation}
for all $Z,W\in \A$ since
\begin{equation}
\label{triple:eq}
ZWZ = 2(Z\circ W)\circ Z - (Z\circ Z)\circ W,
\end{equation}
where $Z\circ W = [(Z+W)^2-Z^2-W^2]/2$.
Let $Z \in G_I(\A)$.  Then $Z^{-1} \in G_I(\A)$ by Example 1. Hence it follows
from~(\ref{ptriple:eq}) and the identities
$I = Z(Z^{-1})^2Z$ and $Z = ZZ^{-1}Z$ that $\rho(Z)^{-1}$ exists and
$\rho(Z^{-1}) = \rho(Z)^{-1}$.  Thus $L(Z)^{-1}$ exists and (\ref{Linv:eq})
holds.  Also,
$L(G_I(\A))$ is a component of $G(\A)$ since $L$ is a homeomorphism.
\hfill $\Box$

{\bf Proof of Theorem \ref{ext:th}.}
Let $L:\A \rightarrow \A$ be an invertible linear isometry with
$L(I) \in G_I(\A)$.  Then by Lemma~\ref{L:lem}, $L(G_I(\A)) = G_I(\A)$
so $L$ is a biholomorphic mapping of $\E = \A_0 \cap G_I(\A)$.
By Example 1, $T(Z) = Z^{-1}$ is a biholomorphic mapping of $\D$ onto
$\E$ and hence
$g = T \circ L \circ T$ is a biholomorphic mapping of $\D$.
Moreover, $L(Z) = U g(Z) U$ for all $Z \in G_I(\A)$ by~(\ref{Linv:eq}).
Therefore, $L$ is a biholomorphic mapping of $\D$.

Now let $h$ be a biholomorphic mapping of $\D$.  Then
$g = T \circ h \circ T$ is a biholomorphic mapping of $\E$.
By hypothesis and the argument given in the proof of Theorem~\ref{LU:th},
$g$ extends to a biholomorphic mapping $\hat{g}$ of $\A_0$.
Put $B = -\hat{g}^{-1}(0)$.  Clearly
$B \in \A_0$ and it follows from Theorem 2.2(ii) and Theorem 2.6
of~\cite{arazy} that
$B^* \in \A$.  Hence
\begin{equation}
\label{TB:eq}
T_B(Z) = (I-BB^*)^{-\frac{1}{2}}(Z+B)(I+B^*Z)^{-1}(I-B^*B)^\frac{1}{2}
\end{equation}
is a biholomorphic mapping of $\A_0$ by~(\ref{triple:eq}),
the comments after Proposition~5 of \cite{LFT} and \cite[(12)]{BSD}.
Thus the proof of Theorem~3 of~\cite{BSD} applies to show that
$\hat{g} = L \circ T_B$,
where $L:\A \rightarrow \A$ is an invertible linear isometry.
By Lemma~\ref{L:lem}, $L^{-1}$ takes $G_I(\A)$ into $G(\A)$ so 
$Z+B$ is invertible whenever $Z \in G_I(\A)$ and
$||Z|| < 1$.  We will show that $B = 0$.  Then $\hat{g} = L$ and
$h = T \circ g \circ T$ so by Lemma~\ref{L:lem}, $h = L_1$ on $\D$, where
$L_1(Z) = VL(Z)V$ and $V=L(I)^*$ is unitary.  Hence $L_1$ is a linear isometry on $\A$
and $L_1(\A) = \A$ since $L_1 = Dh(2I)$.  Moreover, $L_1(I) \in G_I(\A)$
since $h(2I) \in G_I(\A)$.

  To show that $B=0$, let $U$ be any unitary operator in $G_I(\A)$ and note
that $U^* \in G_I(\A)$ since $U^* = T(U)$.  By what we have shown above,
$(-\lambda U^*+B)^{-1}$ exists for all $|\lambda| < 1$ so $|UB|_{\sigma} = 0$,
where $|~~|_{\sigma}$ denotes the spectral radius.
The straightforward extension of \cite[Prop. 2a]{BSD} to $\A$
shows that $T_{B^*}(\lambda I)$ is a unitary operator in $G_I(\A)$ for
all $|\lambda| = 1$.
Thus $f(\lambda) = T_{B^*}(\lambda I) B$ is an operator-valued
function which is holomorphic on a neighborhood of the closed unit disc with
$|f(\lambda)|_{\sigma} = 0$ for all $|\lambda| = 1$ so $|f(0)|_{\sigma} = 0$
by a maximum principle of Vesentini~\cite{vesen}.  Therefore,
$|B^*B|_{\sigma} = 0$ so $B = 0$, as required. \hfill $\Box$

{\bf Proof of Proposition \ref{RSJ:prop}.}
Put $\E = \A_0 \cap G_I(\A)$ and let $h:\E \rightarrow \X$ be a holomorphic
function with $||h(Z)|| \leq 1$ for all $Z \in \E$.  Given $0<t<1/2$, put
$r = 1-t$.  Let $W\in\A_0$ and suppose $\sigma(W)$ is a finite set.  Then
$f(\lambda) = h(tI+\lambda W)$ is a vector-valued function with unit bound
which is holomorphic at all but a finite number of points in the disc
$|\lambda| < r$ so $f$ has a holomorphic extension to this disc
by~\cite[Th. 3.13.3]{HP}.  By the Cauchy estimates~\cite[(3.11.3)]{HP},
$|f^{(n)}(0)| \leq n! r^{-n}$ so
\begin{equation}
\label{CE:eq}
||\hat{D}^n h(tI)(W)|| \leq n! r^{-n},
\end{equation}
where $\hat{D}^n h(Z)$ denotes the homogeneous polynomial associated
with the $n$th order Fr\'{e}chet derivative at $Z$.

By the spectral theorem, each normal operator in $\A$ is a limit in
the operator norm of a linear combination $W$ of orthogonal projections
in $\A$ which sum to $I$.  Hence~(\ref{CE:eq}) holds for all normal operators
$W \in \A_0$.  By the maximum principle for unitary operators~\cite[Th. 9]{BSD},
inequality~(\ref{CE:eq}) holds for all $W\in\A_0$.

Define
\[
h_t(Z) = \sum_{n=0}^{\infty} \frac{1}{n!}\hat{D}^nh(tI)(Z-tI)
\]
and put $\D_t =\{Z\in\A:~||Z-tI|| < r\}$.  Then $h_t$ is holomorphic in
$\D_t$ since by~(\ref{CE:eq}) the above series converges uniformly
on each ball about $tI$ with radius less than $r$. 
Also, $\E$ is connected and $h_t(Z) = h(Z)$ for all $Z \in \A$ with
$||Z-tI|| < t$
by Taylor's theorem~\cite[Th. 3.17.1]{HP}.
Since $\{\D_t:~0<t<1/2\}$ is a monotone family of domains with union $\A_0$,
it follows from the identity
theorem~\cite[Th. 3.16.4]{HP} that there exists a holomorphic extension
$\hat{h}$ of $h$ to $\A_0$.   \hfill $\Box$

In the remainder of this section, we give some types of circular domains which
are homogeneous under the group of linear automorphisms.
Our first result reduces to \cite[Example 4]{LFT} when $H = \C$.

\begin{prop}
\label{prod:prop}
Let $\A = \Lin{H,K\times H}$ and put
\[
\D = \{Z\in\A:~Z_1^*Z_1 < Z_2^*Z_2,~Z_2^{-1}~{\rm exists}\}.
\]
Then $\D$ is symmetric and the group of linear biholomorphic mappings of $\D$ acts transitively
on $\D$.  Moreover, $\D$ is holomorphically equivalent to the product
domain $\D' = \Lin{H,K}_0 \times G\Lin{H}$, where $G\Lin{H}$ is the set
of all invertible operators in $\Lin{H}$.
\end{prop}
{\bf Proof.} Note that
$\D = \{Z\in\E:~Z^*JZ < 0\}$, where
$\E = \{Z\in\A:~Z_2^{-1}~{\rm exists}\}$
and $J = \opmatrix{I}{0}{0}{-I}$.
Given $W\in\D$, put $B = W_1W_2^{-1}$.  Since $||B|| < 1$,
we may define $L(Z) = MZR$, where
\[
M = \opmatrix{(I-BB^*)^{-\frac{1}{2}}} {B(I-B^*B)^{-\frac{1}{2}}}
             {(I-B^*B)^{-\frac{1}{2}}B^*} {(I-B^*B)^{-\frac{1}{2}}}
\]
and $R = (I-B^*B)^{1/2}W_2$.
(Here M is the coefficient matrix of~(\ref{TB:eq}).)
Clearly $L\col{0}{I} = W$.  Moreover, $M^*JM = J$ so
$L(\D) \subseteq \D$ since $L(Z)^*JL(Z) = R^*Z^*JZR$ and
$(I-B^*B)^{-1/2}(B^*Z_1Z_2^{-1} + I)Z_2R$ is invertible for all
$Z\in\D$.  Since $M^{-1}$ is obtained from $M$ by replacement of $B$ by
$-B$ in the expression for $M$, a similar argument shows that
$L^{-1}(\D) \subseteq \D$.  Therefore, $L$ is a biholomorphic
mapping of $\D$.

It is easy to verify that $h(Z) = (Z_1Z_2^{-1},Z_2^{-1})$ is a
biholomorphic mapping of $\D$ onto $\D'$ and hence $\D$ is symmetric since
each of the factors of $\D'$ is symmetric.  \hfill $\Box$

\begin{prop}
\label{prop:hypdom}
Let $J \in \Lin{H}$ with $J^* = J$ and suppose $1$ and $-1$ are eigenvalues
for $J$.  Then the group $G$ of linear biholomorphic mappings of
$\D = \{z \in H:~(Jz,z)<0\}$
acts transitively on $\D$.
\end{prop}

{\bf Proof.} By hypothesis,
$H = \C e \bigoplus K \bigoplus \C f$,
where $e$ and $f$ are unit eigenvectors for $1$ and $-1$, respectively, and
$K = \{e,f\}^{\bot}$.  Thus if $z \in H$, we may write
\[
z = \colthree{\alpha}{w}{\beta}, ~~J=\opmthree{1}{0}{0}{0}{B}{0}{0}{0}{-1},
\]
\[
\D = \{z\in H:~|\alpha|^2 + (Bw,w) < |\beta|^2\}.
\]
Suppose $z_1 \in \D$ and let $\alpha_1$, $w_1$ and $\beta_1$ be the
coordinates of $z_1$.  We will obtain an $L \in G$ with $Lf = z_1$.

First suppose $A \equiv |\alpha_1| + |\beta_1| > 0$ and put
\[
a = |\alpha_1| + \frac{1}{2A}(Bw_1,w_1),~~~
b = |\beta_1| - \frac{1}{2A}(Bw_1,w_1).
\]
Clearly $a + b = A$ and $A(b - a) > 0$ so $b^2 -a^2 > 0$.
Given $w \in K$, let $L_w:H \rightarrow H$ be the linear map defined by the operator matrix
\[
L_w = \opmthree{1-\frac{1}{2}(Bw,w)}{-w^*B}{-\frac{1}{2}(Bw,w)}
                {w}{I}{w}
                {\frac{1}{2}(Bw,w)}{w^*B}{1 + \frac{1}{2}(Bw,w)}.
\]
Then $L_w^{-1} = L_{-w}$ and $L_w^*JL_w = J$ so $L_w \in G$.  Take
$w = w_1/A$.  Then by a short computation,
\[
L_w \colthree{a}{0}{b} = \colthree{|\alpha_1|}{w_1}{|\beta_1|}.
\]
Note that
$L_1 = \opmthree{b}{0}{a}
                {0}{\sqrt{b^2-a^2}}{0}
                {a}{0}{b}$
satisfies $L_1f = \colthree{a}{0}{b}$ and $L_1$
is in $G$ since $L_1^*JL_1 = (b^2-a^2)J$ and $L_1^{-1}$ exists.
Also, there is an $L_2 \in G$
(given by a diagonal matrix) such that
$L_2 \colthree{|\alpha_1|}{w_1}{|\beta_1|} =
      \colthree{\alpha_1}{w_1}{\beta_1}.$
Hence $L_2 L_w L_1 f = z_1$.

Now suppose $\alpha_1 = \beta_1 = 0$.  Then
$L_{w_1} z_1 = \colthree{-(Bw_1,w_1)}{w_1}{(Bw_1,w_1)}$
and $(Bw_1,w_1) \not = 0$ since $z_1 \in \D$.  By what we have shown,
there is an $L \in G$ with $Lf = L_{w_1} z_1$ so $z_1 = L_{w_1}^{-1} Lf$.
  \hfill $\Box$


\noindent
Mathematics Department
\newline
University of Kentucky
\newline
Lexington, Kentucky  40506

\begin{thebibliography}{99}

\bibitem{alpay}
D. Alpay, A. Dijksma, J. van der Ploeg, H.S.V. de Snoo,
{\em Holomorphic operators between Krein spaces and the number of squares of
associated kernels}, in Operator Theory and Complex Analysis,
T. Ando and I. Gohberg, eds., Operator Theory: Advances and Applications,
Vol. 59, Birkh\"{a}user Verlag, Basel, 1992, pp. 11--29.
\vspace{.05in}

\bibitem{arazy}
J. Arazy and B. Solel, {\em Isometries of non-self-adjoint operator algebras},
J. Funct. Anal. {\bf 90}(1990), 284--305.
\vspace{.05in}

\bibitem{BGR}
J. A. Ball, I. Gohberg and L. Rodman, {\em Interpolation of Rational Matrix
Functions}, Operator Theory: Advances and Applications, Vol. 45,
Birkh\"{a}user Verlag, Basel, 1990.
\vspace{.05in}

\bibitem{bogdan}
E. Ligocka and J. Siciak, {\em Weak analytic continuation}, Bull. Acad. Polon.
Sci. S\'{e}r. Sci. Math. Astronom. Phys. {\bf 20}(1972), p. 461--466.
\vspace{.05in}


\bibitem{RGE}
R. G. Douglas, {\em On majorization, factorization, and range inclusion of
operators on Hilbert space}, Proc. Amer. Math. Soc. {\bf 17}(1966), 413--415.
\vspace{.05in}

\bibitem{BAT}
R. G. Douglas, {\em Banach Algebra Techniques in Operator Theory},
Academic Press, New York, 1972.
\vspace{.05in}

\bibitem{franz}
T. Franzoni, {\em The group of holomorphic automorphisms in certain
$J^*$-algebras}, Ann. Mat. Pura Appl. {\bf 127}(1981), 51--66.
\vspace{.05in}

\bibitem{GGK}
I. Gohberg, S. Goldberg and M. Kaashoek, {\em Classes of Linear Operators
Vol. I}, Operator Theory: Advances and Applications, Vol. 49,
Birkh\"{a}user Verlag, Basel, 1990.
\vspace{.05in}

\bibitem{SL}
L. A. Harris, {\em Schwarz's lemma in normed linear spaces}, Proc. Nat. Acad.
Sci. U.S.A. {\bf 62}(1969), 1014--1017.
\vspace{.05in}

\bibitem{BSD}
L. A. Harris, {\em Bounded symmetric homogeneous domains in infinite
dimensional spaces}, in Infinite Dimensional Holomorphy, T. L. Hayden and
T. J. Suffridge, eds., Lecture Notes in Mathematics, Vol 364,
p. 13--40, Springer Verlag, Berlin, 1973.
\vspace{.05in}

\bibitem{OSD}
L. A. Harris, {\em Operator Siegel domains}, Proc. Royal Soc. Edinburgh
Sect. {\bf 79A}(1977), 137--156.
\vspace{.05in}

\bibitem{SPS}
L. A. Harris, {\em Schwarz-Pick systems of pseudometrics for domains in normed
linear spaces}, in Advances in Holomorphy, J. A. Barroso, ed.,
North-Holland, Amsterdam, 1979, pp. 345--406.
\vspace{.05in}

\bibitem{GCA}
L. A. Harris, {\em A generalization of $C^*$-algebras}, Proc. London
Math. Soc. {\bf 42}(1981), 331--361.
\vspace{.05in}

\bibitem{LFT}
L. A. Harris, {\em Linear fractional transformations of circular domains in
operator spaces}, Indiana Univ. Math. J. {\bf 41}(1992), 125--147.
\vspace{.05in}

\bibitem{fact}
L. A. Harris, {\em Factorizations of operator matrices}, Linear Algebra
Appl., to appear.
\vspace{.05in}

\bibitem{AFT}
E. Hille, {\em Analytic Function Theory}, Vol. I, Blaisdell, New York, 1959.
\vspace{.05in}

\bibitem{HP}
E. Hille and R. S. Phillips, {\em Functional Analysis and Semi-Groups},
Amer. Math. Soc. Colloq. Publ., Vol. 31, AMS, Providence, 1957.
\vspace{.05in}

\bibitem{hua1}
Loo-Keng Hua, {\em Geometries of matrices. II. Study of involutions in the
geometry of symmetric matrices},
Trans. Amer. Math. Soc. {\bf 61}(1947), 193--228.
\vspace{.05in}

\bibitem{hua2}
Loo-Keng Hua, {\em  Geometries of matrices. III. Fundamental theorems in the
geometries of symmetric matrices},
Trans. Amer. Math. Soc. {\bf 61}(1947), 229--255.
\vspace{.05in}

\bibitem{isidro}
J. M. Isidro and L. L. Stach\'{o}, {\em Holomorphic Automorphism Groups in Banach 
Spaces: An Elementary Introduction},
Math. Studies Vol. 105, North-Holland, Amsterdam, 1985.
\vspace{.05in}

\bibitem{pota}
V. P. Potapov, {\em Linear fractional transformations of matrices}, in
Studies in the Theory of Operators and Their Applications (Russian),
``Naukova Dumka," Kiev, 1979, pp. 75-97; English transl. in
Amer. Math. Soc. Transl. {\bf 138}(1988), 21--35.
\vspace{.05in}

\bibitem{RHD}
R. C. Penney, {\em The structure of rational homogeneous domains in }$\C^n$,
Ann. of Math. {\bf 126}(1987), 389--414.
\vspace{.05in}

\bibitem{ramis}
J.-P. Ramis, {\em Sous-ensembles analytiques d'une vari\'{e}t\'{e} banachique complexe},
Ergebnisse der Math. No. 53, Springer-Verlag, Berlin, 1970.
\vspace{.05in}

\bibitem{MRD}
A. G. Sergeev, {\em On matrix Rinehardt and circled domains}, in Several
Complex Variables: Proceedings of the Mittag-Leffler Institute, 1987--1988,
J. E. Fornaess, ed., Mathematical Notes 38, Princeton Univ. Press,
Princeton, 1993, pp. 573--586.
\vspace{.05in}

\bibitem{upmeier}
H. Upmeier, {\em Symmetric Banach Manifolds and Jordan $C^*$-Algebras},
Math. Studies Vol. 104, North-Holland, Amsterdam, 1985.
\vspace{.05in}

\bibitem{vesen}
E. Vesentini, {\em On the subharmonicity of the spectral radius},
Boll. Un. Mat. Ital. {\bf 4}(1968), 427--429.
\vspace{.05in}

\bibitem{vigue}
J.-P. Vigu\'e, {\em Le groupe des automorphismes analytiques d'un domaine
born\'{e} d'un espace de Banach complexe.  Application aux domaines born\'{e}s
symmetriques}, Ann. Sci. Ecole Norm. Sup. {\bf 9}(1976), 203--282.
\vspace{.05in}

\bibitem{wink}
J. Winkelmann, {\em On automorphisms of complements of analytic subsets in }$\C^n$,
Math. Z. {\bf 204}(1990), 117--127.
\vspace{.05in}

\end{thebibliography}
\end{document}